\documentstyle[12pt]{amsart}
\newcommand{\der}{\partial}
\newtheorem{thm}{Theorem}[section]
\newtheorem{crl}[thm]{Corollary}
\newtheorem{prp}[thm]{Proposition}
\newtheorem{lm}[thm]{Lemma}

\begin{document}
\title{Representations of $n$-Lie algebras }
\author{A.S. Dzhumadil'daev}
\address{Institut des Hautes \'Etudes Scientifiques, Bures-sur-Yvette
\newline
Institute of Mathematics, Academy of Sciences
of Kazakhstan, {\hbox{Almaty}}}
\email{askar@@math.kz}

\maketitle

\begin{abstract}
Let $V_n=<e_1,\ldots,e_{n+1}>$ be a vector products $n$-Lie algebra
with $n$-Lie commutator $[e_1,\ldots,\hat{e_i},\ldots,e_{n+1}]=(-1)^ie_i$
over the field of complex numbers.
Any finite-dimensional $n$-Lie $V_n$-module
is completely reducible. Any finite-dimensional irreducible
$n$-Lie $V_n$-module is isomorphic to
a $n$-Lie extension of $so_{n+1}$-module with highest weight $t\pi_1$
for some nonnegative integer $t.$
\end{abstract}

\section{Introduction}
An algebra $A=(A,[\;,\;])$ with a $n$-multiplication
$[\;,\;]:\wedge^nA\rightarrow A, (a_1,\ldots,a_n)\mapsto [a_1,\ldots,a_n]$
is called {it $n$-Lie,} if
$$[a_1,\ldots,a_{n-1},[a_n,\ldots,a_{2n-1}]]=
\sum_{i=n}^{2n-1}(-1)^i [[a_1,\ldots,a_{n-1},a_i],
a_n,\ldots,\hat{a_i},\ldots,a_{2n-1}],$$
for any $a_1,\ldots,a_{2n-1}\in A$ \cite{Filipov}.

To any $n$-Lie algebra one can associate
a Lie algebra $L(A)=\wedge^{n-1}A,$ called  {\it basic} Lie algebra,
with a multiplication given by
$$[a_1\wedge \cdots\wedge a_{n-1}, b_1\wedge\cdots\wedge b_{n-1}]=
\sum_{i=1}^{n-1} (-1)^{i+1}[[a_1,\ldots,a_{n-1},b_i],b_1,\ldots,\hat{b_i},
\ldots,b_{n-1}],$$ or by
$$[a_1\wedge \cdots\wedge a_{n-1}, b_1\wedge\cdots\wedge b_{n-1}]=
\sum_{i=1}^{n-1} (-1)^{i+n}[a_1,\ldots,\hat{a_i},\ldots,a_n,
[a_i,b_1,\ldots,b_{n-1}]],$$
where $\hat{b_i}$ means that the element $b_i$ is omitted.

{\bf Example 1.} Let $A=K[x_1,\ldots,x_n]$
under Jacobian map $(a_1,\ldots,a_n)\mapsto det\,(\der_i(a_j)).$
Then $A$ is $n$-Lie \cite{Filipov}, \cite{Filipov1}
and its basic algebra is isomorphic to divergenceless
vector fields algebra $S_{n-1}$ \cite{Dzhu}.

{\bf Example 2.} Let $V_n$ be $(n+1)$-dimensional vector space with a basis
$\{e_1,\ldots,e_{n+1}\}.$ Then $V_n$ under a $n$-Lie multiplication
$[e_1,\ldots,\hat{e_i},\ldots,e_{n+1}]=(-1)^ie_i$
can be endowed by a structure of $n$-Lie algebra. This algebra is called
{\it vector products $n$-Lie algebra}. For $n=2$ we obtain well known
vector products algebra on ${K}^3.$
>From results of \cite{Filipov} it follows that $L(V_n)\cong so_{n+1}.$

Let $U(L(A))$ be an universal enveloping algebra of the Lie algebra $L(A).$
Let $Q(A)$ be an ideal of $U(L(A))$ generated by elements
$$X_{a_1,\ldots,a_{2n-2}}=[a_1,\ldots, a_n]\wedge a_{n+1}\wedge
\cdots \wedge a_{2n-2}$$
$$-\sum_{i=n+1}^{2n-2}(-1)^{i+n}
(a_1\wedge \cdots \hat{a_i}\wedge\cdots\wedge a_n) (a_i\wedge
a_{n+1} \cdots a_{2n-2}),$$ where $a_1,\ldots,a_{2n-2}\in A.$ Let
$U(A)$ be a factor-algebra of $U(L(A))$ by the ideal $Q(A).$

A module $M$ of the Lie algebra $L(A),$ or $U(L(A))$-module, is called {\it
module of $n$-Lie algebra $A$,} if $Q(A)$ acts on $M$ in a trivial way.
Such $n$-Lie modules have the following property:
if $A$ is any $n$-Lie algebra and $M$ is any $n$-Lie module
over $A$, then their semi-direct sum $A+M$ is once again $n$-Lie algebra.
This is the main reason why we consider $L(A)$-modules with the
condition $Q(A)M=0$ as $n$-Lie modules.

One can expect that the $n$-Lie algebra $V_n$  plays in theory
of $n$-Lie algebras a role like $sl_2$ in theory of Lie algebras.
The aim of our paper is to describe all
finite-dimensional representations of vector
products $n$-Lie algebra over the filed of complex numbers.

Let $\pi_1,\ldots,\pi_{[n+1/2]}$ be the fundamental weights for 
$so_{n+1}.$ Recall that $so_4\cong sl_2\oplus sl_2$ and any
irreducible $so_4$-module can be realized as $sl_2\oplus
sl_2$-module $M_{t,r}=M_t\otimes M_r,$ where $M_t$ denotes
$(t+1)$-dimensional irreducible $sl_2$-module with highest weight
$t.$

The main result of our paper is the following

\begin{thm} \label{main} $K={\bf C}.$

{\rm i)} Any finite-dimensional $n$-Lie representation of $V_n, n\ge 2,$
is completely reducible.

{\rm ii)} Let $M_{t,r}$ be an irreducible $so_4$-module with
highest weight $(t,r).$ Then $M_{t,r}$ can be prolonged to $3$-Lie
module over  $V_3,$ if and only if $t=r.$

{\rm iii)} Let $M$ be an irreducible module of Lie algebra
$so_{n+1},$ $n>3,$ with highest weight $\alpha.$  Then $M$ can be
prolonged to $n$-Lie module of $V_n,$ if and only if $\alpha$ has
a form $t\pi_1,$ for some nonnegative integer $t.$
\end{thm}

So, we obtain complete description of finite-dimensional $n$-Lie
$V_n$-modules over ${\bf C}.$ Our result shows that any
irreducible $n$-Lie representation of $V_n$ is ruled by some
nonnegative integer $t$ as in the Lie case $V_2\cong sl_2.$ Call
$t$ mentioned in theorem~\ref{main} {\it $n$-Lie highest weight}.

\begin{crl} $(K={\bf C}, n>2)$
The dimension of any irreducible $n$-Lie $V_n$-module
with highest weight $t$ is equal to $\frac{n+2t-1}{n+t-1}{n+t-1\choose t}.$
\end{crl}

For example, the  dimension of any irreducible
$V_3$-module with highest weight $t$ is equal to $(t+1)^2.$

{\bf Remark.} If $n=3$ and if we allow considering
infinite-dimensional modules, then studying of $V_3$-representations
can be reduced to the problem on describing of $gl_{\lambda}$-modules.
A definition of complex size matrices algebra $gl_{\lambda}$ see
\cite{Dixmier}, \cite{Feigin}.  One can prove that $U(V_3)$ has a
subalgebra isomorphic to $gl_{\lambda}\otimes gl_{\lambda}.$

\section{$n$-Lie modules }

Let $A=(A,\omega)$ be a $n$-Lie algebra with multiplication
$\omega.$ This means that $\omega\in C^n(A,A)$ is a skew-symmetric
polylinear map $\omega: \wedge^nA\rightarrow A$ with $n$ arguments
and satisfies the following $n$-Lie identity (Leibniz rule)
$$\omega(a_1,\ldots,a_{n-1},\omega(a_n,\ldots,a_{2n-1}))=$$
$$\sum_{i=n}^{2n-1}\omega(a_n,\ldots,a_{i-1},
\omega(a_1,\ldots,a_{n-1},a_i),a_{i+1},\ldots,a_{2n-1}),$$
where $a_1,\ldots,a_{2n-1}$ are any elements of $A$ \cite{Filipov1}
The Leibniz rule can be written as follows
$$\omega(a_1,\ldots,a_{n-1},\omega(a_n,\ldots,a_{2n-1}))=$$
$$\sum_{i=n}^{2n-1}(-1)^{n+i}
\omega(\omega(a_1,\ldots,a_{n-1},a_i),a_n,\ldots,\hat{a_i},\ldots,a_{2n-1}),$$
for any $a_1,\ldots,a_{2n-1}\in A.$ Here the notation $\hat{a_i}$
means that $a_i$ is omitted.

Let $End\,A$ be a space of linear maps $A\rightarrow A.$
Recall that an operator $D\in End\,A$ is called {\it derivation,} if
$$D(\omega(a_1,\ldots,a_n))=\sum_{i=1}^n\omega(a_1,\ldots,a_{i-1},D(a_i),a_{i+1},\ldots,a_n),$$
for any $a_1,\ldots,a_n\in A.$
Let $Der\,A$ be a space of derivations of $A.$
According $n$-Lie identity for any $n-1$ elements $a_1,\ldots,a_{n-1}\in A$ one can
correspond adjoint derivation
 $ad\,\{a_1,\ldots,a_{n-1}\}\in Der \,A$ by the rule
$$ad\,\{a_1,\ldots,a_{n-1}\} a=\omega(a_1,\ldots,a_n).$$
Call a derivation $D\in Der\,A$ {\it interior,} if there exists some $a_1,\ldots,a_{n-1}\in A,$
such that $D=ad\,\{a_1,\ldots,a_{n-1}\}.$
Denote by $Int\,A$ a space of interior derivations of $A.$ Then $Der\,A$ is a Lie algebra under
commutator $(D_1,D_2)\mapsto [D_1,D_2]:= D_1D_2-D_2D_1$ and $Int\,A$ is its Lie ideal.
Notice that the commutator of $Int\,A$ can be given by
$$[ad\,\{a_1,\ldots,a_{n-1}\}, ad\,\{b_1,\ldots,b_{n-1}\}]=ad\,\{c_1,\ldots,c_{n-1}\},$$
where
$$\{c_1,\ldots,c_{n-1}\}=\sum_{i=1}^{n-1}(-1)^{i+1}
\omega(\omega(a_1,\ldots,a_{n-1},b_i),
b_1,\ldots,\hat{b_{i}},\ldots,b_{n-1}).$$
By skew-symmetric condition $\{c_1,\ldots,c_{n-1}\}$ can be defined also by
$$\{c_1,\ldots,c_{n-1}\}=-\sum_{i=1}^{n-1}(-1)^i\omega(
\omega(b_1,\ldots,b_{n-1},a_i), a_1,\ldots,\hat{a_i},\ldots,a_{n-1}).$$
Moreover, the $(n-1)$-th exterior power $L(A)=\wedge^{n-1}A$ under commutator
$$[a_1\wedge a_2 \wedge \cdots \wedge a_{n-1},b_1\wedge b_2\wedge \cdots\wedge b_{n-1}]=c_1\wedge c_2\wedge \cdots c_{n-1}$$
can be endowed by structure of Lie algebra
(write $\{a_1,\ldots,a_{n-1}\}$
instead of $a_1\wedge a_2 \wedge \cdots a_{n-1}$ ).
The adjoint map
$$\wedge^{n-1}A\rightarrow Int\,A, \quad a_1\wedge\cdots\wedge a_{n-1}\mapsto ad\,\{a_1,\ldots,a_{n-1}\},$$ gives us
a homomorphism of Lie algebras.
Following L. Takhtajan call $L(A)$ a {\it basic} Lie algebra of $A.$

A vector space $M$ is called  {\it $A$-module,} or {\it 
$n$-Lie module over $A,$} if there are given $n$ polylinear maps
$$\omega_i : A\otimes \cdots \otimes A\otimes \mathop{M}\limits_i \otimes A\cdots \otimes A\rightarrow M, i=1,\ldots,n,$$
such that
\begin{itemize}

\item
$\omega_i(a_1,\ldots,a_{i-1},m,a_{i+1},\ldots,a_n)$ is
skew-symmetric by all $a$-type arguments
\item $$\omega_i(a_1,\ldots,a_{i-1},m,a_{i+1},a_{i+2},\ldots,a_n)=-
\omega_{i+1}(a_1,\ldots,a_{i-1},a_{i+1},m,a_{i+2},\ldots,a_n),$$
for any $i=1,\ldots,n-1,$ and

\item $$\omega_{n}(a_1,\ldots,a_{n-1},\omega_{n}(a_n,\ldots,a_{2n-2},m))=$$
$$\sum_{i=n}^{2n-2}\omega_{n}
(a_n,\ldots,a_{i-1},\omega(a_1,\ldots,a_{n-1},a_i),a_{i+1},\ldots,a_{2n-2},m)$$
$$+\omega_n(a_n,\ldots,a_{2n-2},\omega_n(a_1,\ldots,a_{n-1},m)),$$
for any $a_1,\ldots,a_{2n-2}\in A$ and $m\in M,$

\item
$$\omega_{n-1}(a_1,\ldots,a_{n-2},m,\omega(a_n,\ldots,a_{2n-1}))=$$
$$\sum_{i=n}^{2n-1}\omega_{i}
(a_n,\ldots,a_{i-1},\omega_{n-1}(a_1,\ldots,a_{n-2},m,a_i),a_{i+1},\ldots,a_{2n-1}),
$$
for any $a_1,\ldots,a_{n-2},a_n,\ldots,a_{2n-1}\in A$ and $m\in M.$
\end{itemize}

If $M$ is an $A$-module, we will say that there are given representation
$(\omega_1,\ldots,\omega_n)$
of $A$ on the vector space $M.$

So, any module of $n$-Lie algebra is an usual module of Lie algebra, if $n=2.$
If $n>2,$ then any module of $n$-Lie algebra $A$ is a Lie module of the
basic Lie algebra $L(A)$ under representation
$\rho:\wedge^{n-1}A\rightarrow End\,M$ defined by
$\rho(a_1\wedge \ldots\wedge a_{n-1})(m)=\omega_n(a_1,\ldots,a_{n-1},m),$ such that
$$\rho(\omega(a_1,\ldots,a_n)\wedge a_{n+1}\wedge
\cdots\wedge a_{2n-2})= $$
\begin{equation} \label{n-liemodule}
\sum_{i=1}^n(-1)^{i+n}\rho(a_1\wedge \ldots\wedge \hat{a_i}\wedge
\ldots\wedge a_{n})
\rho(a_i\wedge a_{n+1}\wedge \cdots \wedge a_{2n-2}),
\end{equation}
for any $a_1,\ldots,a_{2n-2}\in A.$ If $M$ is a Lie module over
Lie algebra $L(A)$ that satisfies the condition
(\ref{n-liemodule}) for any $a_1,\ldots,a_{2n-2}\in A,$ then we
will say that Lie module structure on $M$ over $L(A)$ can be {\it
prolonged} to a $n$-Lie module structure over $n$-Lie algebra $A,$
or shortly that Lie module $M$ can be  prolonged to $n$-Lie
module.

{\bf Example.} For any $n$-Lie algebra $A$ its adjoint module,
i.e., a module with vector space $A$ and the action
$(a_1\wedge \ldots \wedge a_{n-1}) b=[a_1,\ldots,a_{n-1},b]$
is $n$-Lie module.

For Lie algebra $L$ denote by $U(L)$ its universal enveloping algebra.
Let $Q(A)$ be an ideal of $U(L(A))$ generated by elements
$$X_{a_1,\ldots,a_{2n-2}}=
\omega(a_1,\ldots,a_n)\wedge a_{n+1}\wedge \cdots\wedge a_{2n-2}$$
$$-\sum_{i=1}^n(-1)^i(a_1\wedge \cdots\wedge \hat{a_i}\cdots\wedge
a_{n}) \, (a_i\wedge a_{n+1}\wedge \cdots\wedge a_{2n-2}).$$ Let
$\bar U(A)$ be factor-algebra of $U(L(A))$ over $Q(A).$

Any Lie module of $L(A)$ can be prolonged to $n$-Lie module, if
and only if it is trivial $Q(A)$-module. In other words, there are
one-to one correspondence between $n$-Lie modules and $\bar
U(A)$-modules. In this sense $\bar U(A)$ can be considered as 
universal enveloping algebra of $n$-Lie algebra $A.$

Let $A$ be a $n$-Lie algebra and $M$ be an  $A$-module. For
simplifying notations we will write
$\omega(a_1,\ldots,a_{i-1},m,a_{i+1},\ldots,a_n)$ instead of
$\omega_i(a_1,\ldots,a_{i-1},m,a_{i+1},\ldots,a_n).$ Then $n$-Lie
module conditions can be written in a unique way
$$\omega(b_1,\ldots,b_{n-1},\omega(b_n,\ldots,b_{2n-1}))=$$
$$\sum_{i=n}^{2n-1}\omega(b_n,\ldots,b_{i-1},\omega(b_1,\ldots,b_{n-1},b_i),
\ldots,b_{i+1},\ldots,b_{2n-1}),$$
for any $b_1,\ldots,b_{2n-1}\in A+M.$

Let $\tilde A=A+M$ be a direct sum of vector spaces $A$ and $M.$
Endow $\tilde A$ by the structure of $n$-multiplication $\tilde\omega$ by
the following rules
$$\tilde\omega(a_1,\ldots,a_n)=\omega(a_1,\ldots,a_n),$$
if $a_1,\ldots,a_n\in A,$
$$\tilde\omega(a_1,\ldots,a_{i-1},m,a_{i+1},\ldots,a_n)=
\omega_i(a_1,\ldots,a_{i-1},m,a_{i+1},\ldots,a_n),$$
if $a_1,\ldots,a_{i-1},a_{i+1},\ldots\in A,$ $m\in M,$ $i=1,\ldots,n,$
and
$$\tilde\omega(b_1,\ldots,b_n)=0,$$
if at least two elements of the set  $\{b_1,\ldots,b_n\}\subset Q$
belongs to $M.$

It is easy to check that $\tilde\omega$ is skew-symmetric and that
$\tilde\omega$ satisfies the Leibniz rule.  So, $\tilde A$ is $n$-Lie, if
$A$ is $n$-Lie and $M$ is $n$-Lie module over $A.$ Call $\tilde A$ {\it semi-direct
sum} or {\it split extension}  of $n$-Lie algebra $A$ by $A$-module $M.$

Suppose that $Q$ is a $n$-Lie algebra with abelian
ideal $M,$ such that $Q/M$ is $n$-Lie. Then $M$ has natural structure of
$n$-Lie module over factor-algebra $Q/M.$
Suppose that $Q$ has subalgebra $A$ isomorphic to $Q/M.$
Then  $Q$ is isomorphic to the semi-direct sum $A+M.$

Let $M$ be a $n$-Lie module over $n$-Lie algebra $A.$ Let $M_1$ be
a subspace of $M,$ such that
$\omega_i(a_1,\ldots,a_{i-1},m,a_{i+1},\ldots,a_{2n-1})\in N,$ for
any $m\in M_1,$ $i=1,\ldots,n,$ and
$a_1,\ldots,\hat{a_i},\ldots,a_{2n-1}\in A.$ In such case we will
say that $M_1$ is {\it $n$-Lie submodule} or just submodule of
$M.$ Any module has trivial submodules $0$ and $M.$ Call $M$ {\it
irreducible,} if any its submodule is trivial. Call $M$ {\it
completely reducible,} if it can be decomposed to a direct sum of
irreducible submodules. Equivalent definition: $M$ is completely
reducible, if for any submodule $N\subseteq M$ one can find
additional submodule $S\subseteq M,$ such that $M\cong N\oplus S.$

\begin{prp} \label{bir}
Let $M$ be a $n$-Lie module over $n$-Lie algebra $A.$ Then any
submodule, any factor-module and dual module of $M$ are also
$n$-Lie modules. If $M_1$ and $M_2$ are $n$-Lie modules over $A,$
then their direct sum $M_1\oplus M_2$ is also $n$-Lie module.
\end{prp}

{\bf Proof.}
Consider $M$ as a Lie module of Lie algebra $L(A).$  Suppose that
$N$ is a Lie submodule of $M.$ Denote by $\rho_M$ and $\rho_N$
representations of $L(A)$ on modules $M$ and $N$ correspondingly.
By our condition $M$ is $n$-Lie module, therefore according
(\ref{n-liemodule}),
$$\rho_M(\omega(a_1,\ldots,a_n)\wedge a_{n+1}\wedge
\cdots\wedge a_{2n-2})= $$
$$\sum_{i=1}^n(-1)^i\rho_M(a_1\wedge \cdots\hat{a_i}\cdots\wedge a_{n})
\rho_M(a_i\wedge a_{n+1}\wedge \cdots\wedge a_{2n-2}),$$
for any $a_1,\ldots,a_{2n-2}\in A.$
Since $\rho_N(a_1\wedge \ldots\wedge a_{n-1})(m)=
\rho_M(a_1\wedge \cdots \wedge a_{n-1})(m),$
for any $a_1,\ldots,a_{n-1}\in A$ and
$m\in N,$ we have
$$\rho_N(\omega(a_1,\ldots,a_n)\wedge a_{n+1}\wedge\cdots\wedge a_{2n-2})= $$
$$\sum_{i=1}^n(-1)^i\rho_N(a_1\wedge \ldots\hat{a_i}\wedge \cdots\wedge a_{n})
\rho_N(a_i\wedge a_{n+1}\wedge \cdots \wedge a_{2n-2}),$$
for any $a_1,\ldots,a_{2n-2}\in A.$
Therefore, $N$ has a structure of $n$-module over $n$-Lie algebra $A.$
By similar reasons $M/N$ is also $n$-Lie module.

It is evident that $M_1\oplus M_2$ under action
$[a_1,\ldots,a_{n-1},m_1+m_2]=[a_1,\ldots,a_{n-1},m_1]+[a_1,\ldots,a_{n-1},m_2]$ is $n$-Lie module.

\begin{crl}\label{eki} Let $M$ be a $n$-Lie module over $n$-Lie algebra $A.$
Then
\begin{itemize}
\item $M$ is irreducible, if and only if $M$ is irreducible as a Lie module over Lie algebra $L(A).$
\item $M$ is completely reducible, if and only if $M$ is completely reducible
as a Lie module over Lie algebra $L(A)$.
\end{itemize}
\end{crl}

{\bf Proof.} The first statement is evident.
Irreducibility of $M$ as a Lie module of $L(A)$
and irreducibility of $M$ as $n$-Lie module are equivalent by definition.

Prove the second one. Let $N$ and $S$ be submodules of $M$, such that
$M\cong N\oplus S.$ Then $N$ and $S$ are Lie submodules of $M$
as modules over the Lie algebra $L(A)$ and $M\cong N\oplus S$
as  Lie modules.
If $M$ is completely reducible as Lie module, then for any Lie submodule
$N\subseteq  M$ one can find another Lie submodule $S$, such that
$M\cong N\oplus S.$ By proposition~\ref{bir} $N$ and $S$ are also
$n$-Lie submodules. So, $M\cong N\oplus S$ as $n$-Lie modules.

\section{Vector products $n$-Lie algebra and its modules}
Let $V_n$ be a vector products $n$-Lie algebra over ${\bf C}.$
It is $(n+1)$-dimensional and the multiplication on a basis
$\{e_1,\ldots,e_n\}$ is  given by
$$[e_1,\ldots,\hat{e_i},\ldots,e_{n+1}]=(-1)^ie_i, \quad i=1,\ldots,n.$$
For example, $V_2$ is the vector products algebra on ${\bf C}^3$ and
as a Lie algebra it is isomorphic to $sl_2.$

Recall that the Lie algebra of skew-symmetric $n\times n$-matrices
$so_n,$ $n\ge 3,$ is semi-simple over $K={\bf C}.$
More exactly, it is simple, if
$n\ne 4$ and has type $B_{[n/2]},$ if $n$ is odd and type $D_{n/2},$ if $n$
is even. If $n=4,$ then $so_4\cong A_1\oplus A_1.$ For $n=3,$ $so_3\cong A_1.$

For $\lambda\in {\bf Q}$ denote by $[\lambda]$ a  maximal integer,
such that $[\lambda]\le \lambda.$ Let $\pi_1,\ldots,\pi_{[n/2]}$
be the fundamental weights of $so_n$ and $M(\alpha)$ be the
irreducible $so_n$-module with highest weight $\alpha.$ Any
highest weight can be characterized by $[n/2]$-typle of
nonnegative integers $\{s_1,\ldots,s_{[n/2]}\},$ namely $\alpha=
\sum_{i=1}^{[n/2]}\alpha_i\pi_i.$ There is another way to describe
highest weights.

Suppose that a sequence of integers or half-integers
$\lambda=\{\lambda_1,\ldots,\lambda_{[n/2]}\}$ satisfies
the following conditions
\begin{itemize}
\item
 $\lambda_1\ge \lambda_2\ge \cdots\ge \lambda_{[n/2]}\ge 0$, if $n$ is odd
and $\lambda_1\ge \lambda_2\ge \cdots\ge \vert\lambda_{n/2}\vert,$ if $n$
is even.
\item
$\alpha_i, i=1,\ldots,[n/2],$ are nonnegative integers, where $\alpha_i
=\lambda_i-\lambda_{i+1}, i=1,\ldots,[n/2]-1$ and $\alpha_{[n/2]}=2\lambda_{[n/2]},$ if $n$ is odd
and $\alpha_{n/2}=\lambda_{n/2-1}+\lambda_{n/2},$ if $n$ is even.
\end{itemize}
Then any irreducible finite-dimensional $so_n$-module with highest weight
$\alpha$  can be restored by a such sequence $\lambda.$

Let $M$ be an irreducible $so_n$-module. For $n\ne 4,$ set
$q(M)=r,$ if its highest weight $\alpha$ satisfies the condition
$\alpha_r\ne 0,$ but  $\alpha_{r'}=0,$ for any
$r'> r.$ For $n=4,$ set $q(M)=1,$ if $so_4$-module is isomorphic to
$M_{t,t},$ for some nonnegative integer $t$ and $q(M)=2,$ if $M\cong M_{t,r},$ for some $t\ne r.$

Let $\alpha$ be a highest weight for $so_n$-module and
$n\ne 4.$ Then  $q(\alpha)=1,$ if and only if $\alpha$ has the form
$k\pi_1$ for some nonnegative integer $k.$

Any finite-dimensional irreducible $sl_2$-module is isomorphic to
$(l+1)$-dimensional irreducible module $M_l$ with highest weight
$l.$ Recall that any highest weight of $sl_2$ can be identified
with a some nonnegative integer. As we mentioned above $so_4\cong
sl_2\oplus sl_2.$ Any irreducible $so_4$-module $M$ can be
characterized by two nonnegative integers $(t,r).$ Namely,
$M\cong M_{t,r}=M_t\otimes M_r,$ where the action of  $a+b$ on
$m+s,$ where $a$ is an element of the first copy of $sl_2$ and $b$
is an element of the second copy of $sl_2$ and $m'\in M_t, m''\in
M_r,$ is given by
$$(a+b)(m'\otimes m'')= a(m')\otimes m''+m'\otimes b(m'').$$
Notice that in this realization to $so_4$-module $M$ with $q(M)=1$
corresponds the $sl_2\oplus sl_2$-module $M_{t,t},$ for $t\ge 0,
t\in {\bf Z}.$

In \cite{Filipov1} is proved that $V_n$ is simple and
any derivation of $V_n$ is interior.
Therefore, $\wedge^{n-1}V_n\cong Int\,V_n.$
More detailed observation of his proof shows that takes place the following

\begin{thm} \label{filipov} For any $n\ge 2,$
$$Der\,V_n\cong Int\,V_n\cong so_{n+1}.$$
In particular,
$$L(V_2)\cong V_2\cong sl_2,$$
$$L(V_3)\cong V_3\wedge V_3\cong sl_2\oplus sl_2,$$
$$L(V_n)\cong \wedge^{n-1}V_n\cong so_{n+1}, n\ge 4.$$
In other words, take place the following isomorphisms of Lie algebras
$$V_2\cong A_1,$$
$$L(V_3)\cong A_1\oplus A_1,$$
$$L(V_n)\cong B_{n/2}, \mbox{\;if\;\;} n\ge 4 \mbox{\;\;is even},$$
$$L(V_n)\cong D_{(n+1)/2}, \mbox{\;\;if\;} n\ge 4 \mbox{\;\;is odd}.$$
\end{thm}

\begin{lm} \label{31jan} The isomorphism of Lie algebras
$L(V_n)\cong so_{n+1}$ can be given by
$$e_1\wedge \cdots\hat{e_i}\wedge\cdots  \hat{e_j}\wedge \cdots\wedge e_{n+1}\mapsto
(-1)^{i+j+n+1}e_{ij},\quad  i<j.$$
\end{lm}

{\bf Proof.} Easy calculations.

\begin{lm}\label{altinai} Let $M$ be $so_{n+1}$-module.
Define quadratic elements $R_{ijsk}$ of $U(so_{n+1})$ by
$$R_{ijsk}=e_{ij} e_{sk}+e_{is} e_{kj}+e_{ik} e_{js}, \quad
1\le i,j,s,k\le n+1.$$ Then $M$ can be prolonged to $n$-Lie
$V_n$-module, if and only if,  $R_{ijsk}m=0,$ for any $m\in M$ and
$1\le i \le n+1, 1\le j<s<k\le n+1, i\not\in \{j,s,k\}.$
\end{lm}

{\bf Proof.}
Below we use the following notation. If $a,b,c,\ldots$ are some vectors,
then $<a,b,c,\ldots>$ denotes their linear span and $\{a,b,c,\ldots \}$
denotes the set of these elements (order of elements are not important) and by
$(a,b,c,\ldots)$ we denote the vector with components $a,b,c,\ldots$ (order of elements is important).

Notice that $X_{a_1,\ldots,a_{2n-2}}$ is skew-symmetric under
arguments  $a_1,\ldots,a_n$ and $a_{n+1},\ldots,a_{2n-2}.$
Therefore, $X_{a_1,\ldots,a_{2n-2}}=0,$ if dimension of the subspace
$<a_1,\ldots,a_n>$ is less than $n$ or dimension of the subspace $<a_{n+1},\ldots,a_{2n-2}>$ is less than $n-2.$

Suppose that $dim\,<a_1,\ldots,a_n>=n.$

Check that $X_{a_1,\ldots,a_{2n-2}}=0,$ if
$V_n\ne <a_1,\ldots,a_{2n-2}>.$ We can assume that $a_1,\ldots,a_{2n-2}$
are basic vectors. Suppose that
$\{a_1,\ldots,a_n\}=\{e_1,\ldots,\hat{e_i},\ldots,e_{n+1}\}$ for some
$i\in\{1,\ldots,n+1\}.$  Since
$V_n$ does not coincide with the subspace
 $<a_1,\ldots,a_{2n-2}>$ and therefore, its dimension is less than $n+1,$
we have $\{a_{n+1},\ldots,a_{2n-2}\}=\{e_1,\ldots,\hat{e_i},\ldots,
\hat{e_j},\ldots,\hat{e_s},\ldots,e_{n+1}\}$ for some $j, s\ne i, j< s.$
Let for simplicity $a_1=e_1,\ldots,a_{i-1}=e_{i-1},
a_{i}=e_{i+1},\ldots,a_n=e_{n+1}$ and
$(a_{n+1},\ldots,a_{2n-2})=(e_1,\ldots,\hat{e_i},\ldots,\hat{e_j},\ldots,\hat{e_s},\ldots,e_{n+1}).$

We have
$$\omega(a_1,\ldots,a_n)=
\omega(e_1,\ldots,\hat{e_i},\ldots,e_n)=(-1)^{i}e_{i}.$$
Further
$$a_r\wedge a_{n+1}\wedge \cdots\wedge a_{2n-2}=0,$$
if $a_r\ne e_i,e_j,e_s.$
Therefore,
$$(a_1\wedge \cdots\hat{a_r}\wedge \cdots\wedge a_{n})
\, (a_r\wedge a_{n+1}\wedge \cdots\wedge a_{2n-2})=0,$$
if $a_r\ne e_i,e_j,e_s.$

Let $f: \wedge^{n-1}V_n\rightarrow so_{n+1}$ be the isomorphism
of Lie algebras constructed in lemma \ref{31jan}. Prolong it to
the isomorphism of universal
enveloping algebras $f:U(\wedge^{n-1}V_n)\rightarrow U(so_{n+1}).$

Thus,
$$f(\omega(a_1,\ldots,a_n)\wedge (a_{n+1}\wedge\cdots\wedge a_{2n-2}))=$$
$$
f((-1)^ie_i\wedge
e_1\wedge\cdots\hat{e_i}\wedge\cdots\hat{e_j}\wedge \cdots
\hat{e_s}\wedge\cdots\wedge e_{n+1})=$$
$$-f(e_1\wedge \cdots\hat{e_j}\wedge\cdots\hat{e_s}\wedge\cdots\wedge e_{n+1})=$$
$$(-1)^{j+s+n}e_{js}.$$

On the other hand
$$\sum_{r=1}^n(-1)^{r+n}f(a_1\wedge \cdots \hat{a_r}\wedge \cdots\wedge a_{n})
\,f(a_r\wedge a_{n+1}\wedge \cdots\wedge a_{2n-2})=$$

$$(-1)^{n+j-1}f(e_1\wedge \cdots\hat{e_i}\cdots \hat{e_j}\cdots
\wedge e_{n+1}) f(e_j\wedge e_1\wedge \cdots\hat{e_i}\cdots\hat{e_j}
\cdots\hat{e_s}\cdots \wedge e_{n+1})$$
$$+(-1)^{n+s-1}f(e_1\wedge \cdots\hat{e_i}\cdots \hat{e_s}\cdots
\wedge e_{n+1}) f(e_s\wedge e_1\wedge \cdots\hat{e_i}\cdots\hat{e_j}
\cdots\hat{e_s}\cdots \wedge e_{n+1})=$$
$$(-1)^{j+s+n+1}e_{ij} e_{is}+(-1)^{j+s+n}e_{is}e_{ij}=$$
$$-(-1)^{j+s+n}[e_{ij},e_{is}]=(-1)^{j+s+n}e_{js}.$$
Therefore, $f(X_{a_1,\ldots,a_{2n-2}})=0,$ and
$X_{a_1,\ldots,a_{2n-2}}=0,$ if the subspace generated by
$a_1,\ldots,a_{2n-2}$ does not coincide with $V_n.$

Now suppose that $V_n$ is generated by elements
$a_1,\ldots,a_{2n-2}.$ As above we can assume that these elements
are basic elements and
$(a_1,\ldots,a_n)=(e_1,\ldots,\hat{e_i},\ldots,e_{n+1})$ and
$(a_{n+1},\ldots,a_{2n-2})=(e_1,\ldots,\hat{e_j},\ldots,\hat{e_s},
\ldots,\hat{e_k},\ldots,e_{n+1})$ for some $1\le i\le n+1, 1\le
j<s<k\le n+1, i\not\in \{j,s,k\}.$ Then
$$\omega(a_1,\ldots,a_n)\wedge a_{n+1}\wedge \cdots \wedge a_{2n-2}=0,$$
since $e_i\in \{a_{n+1},\ldots,a_{2n-2}\}.$ Calculations as above show that
$$\sum_{r=1}^n(-1)^{r+n}f(a_1\wedge \cdots\hat{a_r}\wedge\cdots\wedge a_{n})
\,f(a_r\wedge a_{n+1}\wedge \cdots\wedge a_{2n-2})=$$
$$\pm R_{ijsk}.$$
So, $f(X_{a_1,\ldots,a_{2n-2}})\in <R_{ijsk}: 1\le i\le n+1, 1\le j,s,k\le n+1>.$

Notice that $R_{ijsk}$ is skew-symmetric by arguments $j,s,k.$
Moreover,
$$R_{ijsk}=-e_{ij}[e_{is},e_{ik}]
-e_{is}[e_{ik},e_{ij}]-e_{ik}[e_{ij},e_{is}],$$
and,
$$R_{ijsk}=0,$$
if $i\in \{j,s,k\}.$

So, $so_{n+1}$-module $M$ can be prolonged to $n$-Lie module, if
and only if $R_{ijsk}m=0,$ for any $m\in M, 1\le i\le n+1, 1\le
i<j<k\le n+1.$ $\square$

\medskip

Denote by
$U_r(L)=<x_1 x_2 \ldots x_{r'}: x_i\in L, r'\le r>$
a  subspace of $U(L)$ generated by products of $r$ elements.
Let $0\subset U_1(L)=L\subset U_2(L)\subset
\cdots$ be an increasing filtration of $U(L)$ and
$S(L)$ be a symmetric algebra of $L.$ Endow $U(L)$ and $S(L)\cong gr\,U(L)$
by a structure of $L$-modules under adjoint action.
Then
$S(L)=\oplus_{r\ge 0} S^r(L),$ where
$S^r(L)\cong U_r(L)/U_{r-1}(L)$ is a symmetric $r$-th power of
adjoint representation. In particular, there is a
natural imbedding $Q_2(A)=Q(A)\cap U_2(L(A))\rightarrow S^2(Ad).$
One can establish the following isomorphisms of $so_{n+1}$-modules
$$Q_2(V_n)\cong
\left\{\begin{array}{ll}
M(\pi_4),&\mbox{if\;\;} n\ge 10\\
M(\pi_4+\pi_5),&\mbox{if\;\;} n=9\\
M(2\pi_4),&\mbox{if\;\;} n=8\\
M(2\pi_4)+M(2\pi_3),&\mbox{if\;\;} n=7\\
M(2\pi_3),&\mbox{if\;\;} n=6\\
M(\pi_2+\pi_3),&\mbox{if\;\;} n=5\\
M(\pi_1),&\mbox{if\;\;} n=4\\
{\bf C},&\mbox{if\;\;} n=3\\
\end{array}\right.
$$

Below we use branching rules for irreducible modules corresponding to
the imbedding $so_{n-1}\subset so_n$ given in \cite{Boerner}.

The proof of theorem \ref{main} is based on the following

\begin{thm} \label{main1} Let $k>1.$

{\rm i)} Let $M=M(\alpha)$ be a finite-dimensional
irreducible $so_{2k+1}$-module with highest
weight $\alpha=\sum_{i=1}^k\alpha_i\pi_i.$ Then $M$ as a
module over Lie subalgebra $so_{2k}$
has a submodule, isomorphic to $M(\bar\alpha),$ where
$\bar\alpha=\sum_{i=1}^k\bar\alpha_i\pi_i$ and
$\bar\alpha_i=\alpha_i,$
$i=1,\ldots,k-1,$ and $\bar\alpha_k=\alpha_{k-1}+\alpha_k.$

{\rm ii)} Let $M=M(\alpha)$ be a finite-dimensional irreducible $so_{2k}$-module
with highest weight $\alpha=\sum_{i=1}^k\alpha_i\pi_i.$ Then $M$ as
a module over Lie subalgebra $so_{2k-1}$ has a submodule, isomorphic
to $M(\bar\alpha),$ where $\bar\alpha=\sum_{i=1}^{k-1}\bar\alpha_i\pi_i$
and $\bar\alpha_i=\alpha_i,$
$i=1,\ldots,k-2,$ $\bar\alpha_{k-1}=\alpha_{k-1}+\alpha_{k}.$
\end{thm}

{\bf Proof.}

i) Take
$$\lambda_k=\alpha_k/2, \;
\lambda_{i}=\sum_{j=i}^{k-1}\alpha_{j}+\alpha_k/2, \;  1\le i\le k.$$
According branching theorem 12.1b \cite{Boerner} any $so_{2k}$-submodule of
$M(\alpha)$ has weight of the form $\bar\alpha,$ such that corresponding
$\bar\lambda$ satisfies the following inequality
$$\lambda_1\ge \vert \bar\lambda_1\vert\ge \lambda_2\ge \cdots \ge
\lambda_{k-1}\ge \vert\bar\lambda_{k-1}\vert \ge \lambda_k\ge
\vert\bar\lambda_k\vert.$$
The $\bar\lambda_j$ are integral or half-integral according to what the
$\lambda_j$ are.
If we take  $\bar\lambda_i:=\lambda_i,$ then such $\bar\lambda$ satisfies these conditions.
Therefore, $M(\alpha)$ has $so_{2k}$-submodule
isomorphic to $M(\bar{\alpha}),$ where
$\bar\alpha=\sum_{i=1}^k\bar\alpha_i\pi_i,$
$\bar\alpha_i=\bar\lambda_i-\bar\lambda_{i+1}=\alpha_i,$ for $i=1,\ldots,k-1,$
and $\bar\alpha_k=\bar\lambda_{k-1}+\bar\lambda_k=\lambda_{k-1}+\lambda_k=\alpha_{k-1}+\alpha_k.$
So, the $so_{2k+1}$-module $M(\alpha)$ as $so_{2k}$-module has a submodule
isomorphic to $M(\bar\alpha),$ where  $\bar\alpha=\sum_{i=1}^{k-1}\alpha_i\pi_i+
(\alpha_{k-1}+\alpha_k)\pi_k.$

ii) We have
$$\alpha_i=\lambda_i-\lambda_{i+1}, 1\le i\le k-1, \alpha_k=\lambda_{k-1}+\lambda_k.$$
By branching theorem 12.1a \cite{Boerner} any $so_{2k-1}$-submodule of
$M(\alpha)$ is isomorphic to $M(\bar\alpha),$ such that corresponding
$\bar\lambda$ satisfies the following inequality
$$\lambda_1\ge \vert \bar\lambda_1\vert\ge \lambda_2\ge \cdots \ge
\lambda_{k-1}\ge \vert\bar\lambda_{k-1}\vert \ge
\vert\lambda_k\vert.$$
The $\bar\lambda_j$ are integral or half-integral according to what the
$\lambda_j$ are. Notice that a sequence $\bar\lambda$ constructed by the following way
satisfies these conditions
$$\bar\lambda_i=\lambda_i, \; 1\le i \le n-1.$$
So, $so_{2k}$-module $M(\alpha)$ as $so_{2k-1}$-module has submodule
$M(\bar\alpha),$ where
$$\bar\alpha_i=\bar\lambda_i-\bar\lambda_{i+1}=\alpha_i, 1\le i \le k-2,$$
$$\bar\alpha_{k-1}=2\bar\lambda_{k-1}=2\lambda_{k-1}=\alpha_{k-1}+\alpha_k.$$

\begin{crl} \label{main2}
Let $n> 4$ and $M$ be irreducible $so_n$-module
such that $q(M)>1.$ Then $M$ as a module over subalgebra
$so_{n-1}\subset so_n$ has a submodule isomorphic to $\bar M,$
such that $q(\bar M)>1.$
\end{crl}

{\bf Proof.} It is easy to see  that for irreducible module $M$ with
highest weight $\alpha,$ the condition $q(M)>1,$ is equivalent to the
condition $\sum_{i>1}\alpha_i>0.$

Let $\bar\alpha$ be highest weight of $so_{n-1},$ defined by
$\bar\alpha_i=\alpha_i, i=1,\ldots,k-1,$ and $\bar\alpha_k=\alpha_{k-1}+\alpha_k,$ if
$n=2k+1,$ and  $\bar\alpha_i=\alpha_i, i=1,\ldots,k-2,$
$\bar\alpha_{k-1}=\alpha_{k-1}+\alpha_{k},$ if $n=2k.$

Notice that
$$\sum_{i>1}\bar\alpha_i= \sum_{i>1}\alpha_i+2\alpha_{k-1}\ge \sum_{i>1}\alpha_i,$$
if $n=2k+1, k>1$ and
$$\sum_{i>1}\bar\alpha_i= \sum_{i>1}\alpha_i,$$
if $n=2k, k>2.$

By theorem \ref{main1} $so_n$-module $M=M(\alpha)$ as $so_{n-1}$-module has a
submodule isomorphic to $\bar M=M(\bar\alpha).$ If $q(M)>1,$ then
$q(\bar M)>1.$ $\square$

\medskip

Notice that $gl_n$ can be realized as a Lie algebra of derivations
of $K[x_1,\ldots,x_n]$ of the form
$\sum_{i,j=1}^n\lambda_{ij}x_i\der_j,$ $\lambda_{ij}\in K.$ Its
subalgebra $so_n$ is generated by elements
$e_{ij}=x_i\der_j-x_j\der_i.$ The set $\{e_{ij}: 1\le i<j\le n\}$
consists of basis of $so_n.$ The multiplication on $so_n$ can be
given by
$$[e_{ij},e_{sk}]=0, \mbox{\;\; if \;\;} \vert \{i,j,s,k\}\vert =4,$$
$$[e_{ij},e_{is}]=-e_{js},
[e_{ij},e_{js}]=e_{is},
[e_{is},e_{js}]=-e_{ij}.$$

\begin{lm} \label{30jan}
Let $M=M_{t,r}$ be an irreducible $so_4$-module. Then $M$ can be
prolonged to $3$-module over $3$-Lie vector products algebra
$V_3,$ if and only if $t=r.$
\end{lm}

{\bf Proof.}  The algebra $so_4$ has the basis $\{e_{ij} : 1\le i<j\le 4\}.$
Take here another basis $\{f_i : 1\le i \le 6\},$ by
$$f_1=(e_{12}+e_{34})/2, f_2=(e_{13}-e_{24})/2, f_3=(e_{14}+e_{23})/2,$$
$$f_4=(-e_{12}+e_{34})/2, f_5=(e_{13}+e_{24})/2, f_6=(-e_{14}+e_{23})/2.$$
Then
$$[f_1,f_2]=-f_3, [f_1,f_3]=f_2, [f_2,f_3]=-f_1,$$
$$[f_4,f_5]=f_6, [f_5,f_6]=f_4, [f_6,f_4]=f_5,$$
$$[f_i,f_j]=, \; i=1,2,3, \; j=4,5,6.$$
We see that
$$e_{12}=f_1-f_4, e_{13}=f_2+f_5, e_{14}=f_3-f_6,$$
$$e_{23}=f_3+f_6, e_{24}=-f_2+f_5, e_{34}=f_1+f_4,$$
and
$$R_{1234}=e_{12}e_{34}-e_{13}e_{24}+e_{14}e_{23}=$$
$$f_1^2-f_4^2+f_2^2-f_5^2+f_3^2-f_6^2=$$
$$C_1-C_2,$$
where
$$C_1=f_1^2+f_2^2+f_3^2, \; C_2=f_4^2+f_5^2+f_6^2,$$
are Casimir elements of subalgebras $<f_1,f_2,f_3>\cong sl_2$
and $<f_4,f_5,f_6>\cong sl_2.$ Well known that any irreducible
finite-dimensional $sl_2$-module is uniquely defined by eigenvalue of the
Casimir operator on this module. Therefore, $M_{t,r}$ is $3$-Lie module,
if and only if $t=r.$

\begin{lm} \label{ush}  Let $n>3.$ Any irreducible
$so_{n+1}$-module $M(t\pi_1)$ can be prolonged to $n$-Lie module
of $V_n.$ Let $M$ be an irreducible $so_{n+1}$-module with
$q(M)>1.$ Then $M$ can not be prolonged to $n$-Lie module over
$n$-Lie algebra $V_n.$
\end{lm}

{\bf Proof.} Let $n> 3.$ Let us consider realization of
$M(t\pi_1)$ as a space of homogeneous polynomials $\sum_{1\le
i_1\le\ldots \le i_t\le n+1} \lambda_{i_1\ldots i_t}x_{i_1}\ldots
x_{i_t}.$

By lemma \ref{altinai} we need to check that
$$R_{ijsk} u=0, \mbox{\;\; for\;\;} u=x_{i_1}\ldots x_{i_t},$$
for any $\{i,j,s,k\}$, such that $1\le i\le n+1,$ $1\le j\le s\le
k\le n+1,$ $i\not\in \{j,s,k\}$ and  $1\le i_1\le i_2\le \cdots
\le i_t\le n+1.$

Let $I=\{i,j,s,k\}.$  Present $u$ in the form $v w,$ where
$v=\prod_{l\in I\cap \{i_1,\ldots,i_t\}}x_l$ and $w=\prod_{l\in
\{i_1,\ldots,i_t\}\setminus I}x_l.$ Notice that
$$R_{ijsk}(vw)=R_{ijsk}(v)w.$$
Therefore it is enough to check that $R_{ijsk}(v)=0,$ for elements
$v\in M(t\pi_1)$ of the form
$$v=x_ix_jx_sx_k,  1\le i\le n+1, 1\le j<s<k\le n+1, i\not\in \{j,s,k\},$$
$$v=x_jx_sx_k, 1\le j\le s\le k\le n+1,$$
$$v=x_ix_sx_k, 1\le i\le n+1, 1\le s\le k\le n+1,$$
$$v=x_j x_k, 1\le j\le k\le n+1,$$
$$v=x_i x_k, 1\le i\le n+1, 1\le k\le n+1,$$
$$v=x_k, 1\le k\le n+1,$$
$$v=x_i, 1\le i\le n+1.$$

Let $i\ne j,s,k.$ Then
$$R_{ijsk}(x_ix_jx_sx_k)=$$
$$e_{ij}(x_ix_jx_s^2-x_ix_jx_k^2)+e_{is}(x_i x_sx_k^2-x_ix_j^2x_s)+
e_{ik}(x_ix_j^2x_k-x_ix_s^2x_k)=$$
$$e_{ij}(x_ix_j)x_s^2-e_{ij}(x_ix_j)x_k^2
+e_{is}(x_ix_s)x_k^2-e_{is}(x_ix_s)x_j^2
+e_{ik}(x_ix_k)x_j^2-e_{ik}(x_ix_k)x_s^2=$$
$$(x_i^2-x_j^2)x_s^2-(x_i^2-x_j^2)x_k^2
+(x_i^2-x_s^2)x_k^2-(x_i^2-x_s^2)x_j^2
+(x_i^2-x_k^2)x_j^2-(x_i^2-x_k^2)x_s^2=0,$$

$$R_{ijsk}(x_jx_sx_k)=$$
$$e_{ij}(x_jx_s^2-x_jx_k^2)+e_{is}( x_sx_k^2-x_j^2x_s)+
e_{ik}(x_j^2x_k-x_s^2x_k)=$$
$$e_{ij}(x_j)x_s^2-e_{ij}(x_j)x_k^2+e_{is}(x_s)x_k^2-e_{is}(x_s)x_j^2
+e_{ik}(x_k)x_j^2-e_{ik}(x_k)x_s^2=$$
$$x_ix_s^2-x_ix_k^2+x_ix_k^2-x_ix_j^2+x_ix_j^2-x_ix_s^2=0,$$

$$R_{ijsk}(x_ix_sx_k)=$$
$$e_{ij}(x_ix_s^2-x_ix_k^2)+e_{is}( -x_ix_jx_s)+e_{ik}(x_ix_jx_k)=$$
$$e_{ij}(x_i)x_s^2-e_{ij}(x_i)x_k^2-e_{is}(x_ix_s)x_j+e_{ik}(x_ix_k)x_j=$$
$$-x_jx_s^2+x_jx_k^2-x_i^2x_j+x_s^2x_j+x_i^2x_j-x_k^2x_j=0,$$

$$R_{ijsk}(x_sx_k)=$$
$$e_{ij}(x_s^2-x_k^2)+e_{is}( -x_jx_s)+e_{ik}(x_jx_k)=-x_j^2+x_j^2=0,$$

$$R_{ijsk}(x_ix_k)=$$
$$e_{ij}(x_ix_s)+e_{is}( -x_ix_j)=-x_jx_s+x_sx_j=0,$$

$$R_{ijsk}(x_k)=e_{ij}(x_s)+e_{is}(-x_j)=0,$$

$$R_{ijsk}(x_i)=0.$$
So, we have checked that $Q(V_n)M(t\pi_1)=0,$ if $n>3.$

Suppose now that $q(M)>1.$
We need to prove that $R_{ijsk}m\ne 0,$ for some
$1\le i\le n+1, 1\le j<s<k\le n+1$ and $m\in M.$

Let us use induction on $n\ge 3.$ If $n=3,$ then by
lemma~\ref{30jan} any irreducible $so_{n+1}$-module $M$ with
$q(M)>1$ can not be prolonged to $n$-Lie module. Suppose that the
statement is true for $n-1\ge 3.$ If $q(M)>1$ for
$so_{n+1}$-module $M,$ then by corollary~\ref{main2} there exists
its $so_{n}$-submodule $\bar M,$ such that $q(\bar M)>1.$ Then by
inductive suggestion  there exists some $R_{ijsk}\in
Q(V_{n-1})\subset U(so_{n})$ and $m\in \bar M,$ such that
$R_{ijsk}m\ne 0.$ Since $m\in \bar M\subseteq M$ and $R_{ijsk}\in
U(so_{n})\subset U(so_{n+1}),$ this means that $R_{ijsk}m\ne 0$ as
elements of $M.$ So, we have proved that our statement for $n.$

{\bf Proof of theorem \ref{main}.}

i) By theorem \ref{filipov} Lie algebra  $\wedge^{n-1}V_n\cong
so_{n+1}$ is semi-simple. Therefore, by Weyl theorem and
proposition \ref{eki}, any any finite-dimensional $n$-Lie
representation of $V_n$ is completely reducible.

ii) and iii) For $n=2$ our statements are evident.
Let $n>2.$ By lemma~\ref{ush} and lemma~\ref{30jan}
$M(t\pi_1), n>3,$ or $M_{t,t}, n=3,$
is $V_n$-module for any nonnegative integer $t$
and any module with $q(M)>1$ can not be $n$-Lie module.

\begin{center}
{\em ACKNOWLEDGEMENTS}
\end{center}

I am grateful to INTAS foundation for support.


\begin{thebibliography}{10}

\bibitem{Boerner} H. Boerner, {\em Group representations,} Berlin 1955.

\bibitem{Dixmier} J. Dixmier, {\em Quotients simples de l'alg\'ebre envelopante
de $sl_2,$} J. Algebra, {\bf 24}(1973), 551-564.

\bibitem{Filipov} V.T. Filipov, {\em $n$-Lie algebras,} Sib. Mat.
J. {\bf 26}(1985), 126-140.

\bibitem{Filipov1} V.T. Filipov, {\em On $n$-Lie algebra of
jacobians,} Sib.Mat.J., {\bf 39}(1998), 660-669.

\bibitem{Feigin} B. L. Feigin, {\em The Lie algebra $gl(\lambda)$ and the cohomologies
of differential operators}, Uspechi Mathem. nauk, {\bf 43}(1988),
No.2, 157-158.

\bibitem{Dzhu} A.S. Dzhumadil'daev, {\em Identities and derivations for
jacobian algebras,} Proc. workshop ``Quantization, deformations, and
new homological and categorical methods in mathematical physics'',
Manchester, 2001.
\end{thebibliography}
\end{document}